\documentclass[a4paper,11pt]{article}

\usepackage{amsmath,amssymb,amsthm,a4wide,color}
\usepackage[utf8]{inputenc}
\usepackage[T1]{fontenc}
\usepackage{authblk,mathrsfs}
\usepackage{algorithm,algpseudocode}
\usepackage{graphicx,subcaption}
\usepackage{bm}

\newcommand{\RR}{\mathbb{R}}
\newcommand{\CC}{\mathbb{C}}
\newcommand{\mrm}{\mathrm}
\newcommand{\mL}{\mathrm{L}}
\newcommand{\mH}{\mathrm{H}}
\newcommand{\mJ}{\mathrm{J}}
\newcommand{\mV}{\mathrm{V}}
\newcommand{\mW}{\mathrm{W}}
\newcommand{\Id}{\mathrm{Id}}
\newcommand{\mD}{\mathrm{D}}
\newcommand{\mT}{\mathrm{T}}

\newcommand{\mR}{\mathrm{R}}
\newcommand{\mB}{\mathrm{B}}
\newcommand{\mA}{\mathrm{A}}
\newcommand{\mQ}{\mathrm{Q}}

\newcommand{\mbX}{\mathbb{X}}

\newcommand{\Vh}{\mathrm{V}_{h}}
\newcommand{\mbVh}{\mathbb{V}_{h}}
\newcommand{\Wh}{\mathrm{W}_{h}}
\newcommand{\mbWh}{\mathbb{W}_{h}}

\newcommand{\bff}{\boldsymbol{f}}

\newcommand{\bx}{\boldsymbol{x}}

\newcommand{\bn}{\boldsymbol{n}}
\newcommand{\bu}{\boldsymbol{u}}
\newcommand{\bv}{\boldsymbol{v}}

\newcommand{\be}{\boldsymbol{e}}

\newcommand{\Bh}{\mathrm{B}}

\newtheorem{defn}{Definition}[section]
\newtheorem{lem}[defn]{Lemma}
\newtheorem{thm}[defn]{Theorem}
\newtheorem{prop}[defn]{Proposition}
\newtheorem{cor}[defn]{Corollary}

\title{\textbf{Substructuring the Hiptmair-Xu}\\
  \textbf{preconditioner for positive Maxwell problems}}
\author[*,$\dagger$]{R.Delville-Atchekzai}
\author[*]{X.Claeys}
\author[$\dagger$]{M.Lecouvez}

\affil[*]{\small Sorbonne Université,
    Université Paris-Diderot SPC,
    CNRS, Laboratoire Jacques-Louis Lions}
\affil[$\dagger$]{CEA-Cesta}

\date{}

\begin{document}

\maketitle

\begin{abstract}
  \noindent 
  Considering positive Maxwell problems, we propose a substructured
  version of the Hiptmair-Xu preconditioner based on a new formula
  that expresses the inverse of Schur systems in terms of the inverse
  matrix of the global volume problem.
\end{abstract}

\section*{Introduction}

Although many peconditioning approaches are now available for
positive definite problems stemming from scalar valued PDEs like
Laplace equation and scalar diffusion, current literature still
offers few techniques for dealing with Maxwell type problems, even in
the symmetric positive definite case.

In the context of unstructured finite element discretization, a major
contribution has been made by Hiptmair and Xu \cite{MR2361899} by
combining the auxiliary space approach \cite{MR1393008} with the
concept of regular decomposition of fields \cite[\S 2.4]{MR2009375},
\cite{MR4143285}. The Hiptmair-Xu method relies on a preconditioner
for scalar Laplace-like problems and, hence, naturally lends itself to
volume domain decomposition based on an overlapping decomposition of
the computational domain, or on multi-grid strategies. This is how it
has been considered in all subsequent contributions
\cite{MR4319101,MR4211733,MR3649422,MR2536904,li2023nodal}.

A substructuring domain decomposition approach called BDDC-deluxe was
also developed by Dohrmann and Widlund
\cite{MR2970732,MR3242973,MR3465088,MR3643544} to deal with positive
curl-curl problems. Although part of the analysis of their
contributions relies on results established by Hiptmair and Xu, the
BDDC algorithm itself appears independent of the strategy advocated in
\cite{MR2361899}. It rests on a dual/primal decomposition of the unknowns
on the skeleton of the partition and on a well chosen averaging operator. 

\quad\\
In the present contribution, we derive another substructuring strategy
for positive curl-curl problems that stems directly from the
Hiptmair-Xu preconditioner. After a technical preamble in Section
\ref{TechnicalPreamble} concerning Moore-Penrose pseudo-inverses, in
Section \ref{GeometricConfig} we describe the geometrical and
discretization setting under consideration. Then in Section
\ref{HarmonicLifting}, we discuss the notion of harmonic lifting and
Theorem \ref{FormulaInverse} establishes an explicit formula for the
inverse of a Schur complement system in terms of the inverse of the
matrix of the global volume problem. This formula seems to be
new. It allows to convert any volume based preconditioner into a
substructuring preconditioner for the associated Shur complement
system. This is the path we follow, applying this idea to the
Hiptmair-Xu method in Section \ref{HXPrecSec} where we also show that
the condition number estimates available for the Hiptmair-Xu
preconditioner can be readily transferred to its substructured
variant, see Lemma \ref{FinalCondEstimate}.

\quad\\
As mentioned above, the Hiptmair-Xu preconditioner needs to be
constructed on top of an already available preconditioner for the
underlying scalar Laplace-like operator. In this respect, for the sake
of clarity, we base our analysis on the one-level Neumann-Neumann
preconditioner. We recall it in detail in Section
\ref{RappelNeumannNeumann} and the overall method is illustrated by
numerical results presented in Section \ref{RappelNeumannNeumann} and
\ref{HXPrecSec}. The one-level Neumann-Neumann preconditioner is
certainly not the best preconditioner for scalar Laplace-like
problems, and adding a coarse space contribution would definitely
improve the performances, but it is easy to introduce. The present
contribution does \textit{not} aim at devising the most efficient
substructuring preconditioner for positive curl-curl problems. Our
goal is simply to show how a preconditioner for the Schur complement
system of a positive curl-curl problem can be deduced from a
preconditioner for the Schur complement system of an associated
Laplace problem.

\paragraph{Notation conventions} In the following when $\mH$
refers to an Hilbert space, then $\mH^{*}$ refers to
its dual i.e. the space of bounded linear forms over $\mH$.
The canonical duality pairing shall be systematically denoted
$\langle \cdot,\cdot\rangle:\mH^{*}\times \mH\to \CC$, and we
shall write $\langle \varphi,v\rangle$ or $ \langle v,\varphi\rangle$
to refer to $\varphi(v)$ for any $\varphi\in \mH^*,v\in\mH$.
Duality pairing shall not involve any complex conjugation.
If $\mV$ is another Hilbert space and $\mA:\mV\to \mH$ is a bounded map,
then $\mA^*:\mH^*\to \mV^*$ shall refer to the dual map defined
by $\langle \mA^*(\varphi), u\rangle := \langle \varphi, \mA(u)\rangle$
for all $\varphi\in \mH^*, u\in \mV$.

\section{Preamble on weighted pseudo-inverses}\label{TechnicalPreamble}

As a preamble, we start by discussing the notion Moore-Penrose
pseudo-inverse (also called generalized inverse) that plays an
important role in the forthcoming analysis. For a comprehensive
presentation see \cite{zbMATH01889793}, see also
\cite{zbMATH03189639,zbMATH03302838}.  Consider two Hilbert spaces
$\mV,\mH$, assuming that the scalar product over $\mH$ is induced by
an operator $\mA:\mH\to \mH^*$ and the associated norm shall be
denoted
\begin{equation*}
  \Vert u\Vert_{\mA}:=\vert\langle \mA u,\overline{u}\rangle\vert^{1/2}
\end{equation*}
Next consider a bounded surjective (but a priori not injective)
operator $\Theta:\mH\to \mV$.  Then the Moore-Penrose pseudo-inverse of
$\Theta$ relative to $\mA$ refers to the linear operator
$\Theta^{\dagger}_{\mA}:\mV\to \mH$ defined, for all $u\in\mV$, by
\begin{equation}\label{DefPseudoInverse}
  \begin{aligned}
    & \Theta\cdot\Theta^{\dagger}_{\mA}(u) = u\quad \text{and}\\
    & \Vert\Theta^{\dagger}_{\mA}(u)\Vert_{\mA}=\inf\{\;\Vert v\Vert_{\mA}:\;v\in\mH,\;\Theta(v) = u\;\}.
  \end{aligned}
\end{equation}
The operator $\Theta^{\dagger}_{\mA}$ definitely depends on the choice
of the scalar product $\mA$. With the previous construction, it is an
injective bounded operator, see e.g. \cite[Lemma
  1]{zbMATH03302838}. We readily check that
$(\Theta^{\dagger}_{\mA}\cdot\Theta)^{2} =
\Theta^{\dagger}_{\mA}\cdot\Theta$
i.e. $\Theta^{\dagger}_{\mA}\cdot\Theta:\mH\to \mH$ is a projector. On
the other hand, since for any $u\in\mH$ we have $u\in
\{\;v\in\mH:\;\Theta(v) = \Theta(u)\;\}$, we obtain that
$\Vert(\Theta^{\dagger}_{\mA}\Theta)u\Vert_{\mA}\leq \Vert
u\Vert_{\mA}\forall u\in \mH$.  Since contractive projectors are
orthogonal, see e.g. \cite[Chap.2 \S 8.1]{MR1192782}, the operator
$\Theta^{\dagger}_{\mA}\Theta$ is self-adjoint in the scalar product
$\mA$ which also writes
\begin{equation}\label{SelfAdjointProjector}
\mA\cdot(\Theta^{\dagger}_{\mA}\Theta) = (\Theta^{\dagger}_{\mA}\Theta)^*\cdot\mA
\end{equation}
Pseudo-inverses are ubiquitous in domain decomposition, at least for
that part of the literature that deals with symmetric positive
definite problem. In this context, the above identity
yields $(\Theta_{\mA}^{\dagger})^*\mA(\Theta_{\mA}^{\dagger})\Theta\mA^{-1}\Theta^* =
(\Theta_{\mA}^{\dagger})^*(\Theta_{\mA}^{\dagger}\Theta)^*\mA\mA^{-1}\Theta^* =
(\Theta\Theta_{\mA}^{\dagger})^*(\Theta\Theta_{\mA}^{\dagger})^* = \Id$, an interesting
identity which we summarize with the following lemma.
\begin{lem}\label{TransfoInverse}\quad\\
  Let $\mV,\mH$ be two Hilbert spaces, $\mA:\mH\to \mH^*$ a bounded self-adjoint
  and coercive operator and $\Theta:\mH\to \mV$ a bounded surjective operator.
  Then $\Theta\mA^{-1}\Theta^*:\mV^*\to \mV$ is a bounded isomorphism, and
  we have $(\Theta\mA^{-1}\Theta^*)^{-1} = (\Theta_{\mA}^{\dagger})^*\mA(\Theta_{\mA}^{\dagger})$.
\end{lem}

\noindent 
Now consider another map $\Phi:\mV\to\mH$ that is
bounded and injective but a priori not surjective. In this situation,
the Moore-Penrose pseudo-inverse $\Phi^{\dagger}_{\mA}:\mH\to
\mV$ is defined, for all $u\in \mH$, by the equations
\begin{equation}
  \begin{aligned}
    & \Phi^{\dagger}_{\mA}\Phi(u) = u\quad \text{and}\\
    & \Vert u-\Phi\Phi^{\dagger}_{\mA}(u)\Vert_{\mA} =
    \inf\{\Vert u-\Phi(v)\Vert_{\mA}: v\in \mV\}.
  \end{aligned}
\end{equation}
Here again $\Phi\Phi^{\dagger}_{\mA}$ is
an $\mA$-orthogonal projector i.e. $\mA\cdot(\Phi\Phi^{\dagger}_{\mA}) =
(\Phi\Phi^{\dagger}_{\mA})^{*}\cdot \mA$. Using this identity
we readily compute
$\Phi^{\dagger}_{\mA}\mA^{-1}(\Phi^{\dagger}_{\mA})^* \Phi^*\mA\Phi =
 \Phi^{\dagger}_{\mA}\mA^{-1}(\Phi\Phi^{\dagger}_{\mA})^*\mA\Phi =
 \Phi^{\dagger}_{\mA}\mA^{-1}\mA \Phi\Phi^{\dagger}_{\mA}\Phi =
 (\Phi^{\dagger}_{\mA}\Phi)^2 = \Id$. We have just proved the following
 lemma. 

\begin{lem}\label{TransfoInverse2}\quad\\
  Let $\mV,\mH$ be two Hilbert spaces, $\mA:\mH\to \mH^*$ a bounded self-adjoint
  and coercive operator and $\Phi:\mV\to \mH$ a bounded injective operator.
  Then $\Phi^*\mA\Phi:\mV\to \mV^*$ is a bounded isomorphism, and
  we have $(\Phi^*\mA\Phi)^{-1} = (\Phi^{\dagger}_{\mA})\mA^{-1}
  (\Phi^{\dagger}_{\mA})^*$.
\end{lem}

\section{Geometric configuration}\label{GeometricConfig}
In the present article we consider a bounded polyhedral computational
domain $\Omega\subset \RR^3$, and a regular triangulation
$\mathcal{T}_{h}(\Omega)$ of $\overline{\Omega} =
\cup_{\tau\in\mathcal{T}_{h}(\Omega)}\overline{\tau}$.  Shape
regularity of this mesh is \textit{not} needed for the subsequent
analysis. The space $\mL^{2}(\Omega)$ refers to square integrable
functions over $\Omega$ and $\mH^{1}(\Omega):=\{
v\in\mL^{2}(\Omega),\; \nabla v\in\mL^{2}(\Omega)\}$.
We consider $\mathbb{P}_1$-Lagrange finite element discretization
\begin{equation*}
  \Vh(\Omega) := \{ v\in \mH^{1}(\Omega),\; v\vert_{\tau}\in
  \mathbb{P}_1(\tau)\;\forall\tau\in\mathcal{T}_h(\Omega) \}
\end{equation*}
Although the results we are going to present can be readily adapted
to higher order $\mathbb{P}_k$-Lagrange discretizations, we
choose to stick to lowest order elements for the sake of clarity.
If $\omega\subset \Omega$ is any open subset that is resolved by the
triangulation i.e. $\overline{\omega} = \cup_{\tau\in
  \mathcal{T}_{h}(\omega)}\overline{\tau}$, where
$\mathcal{T}_{h}(\omega)\subset \mathcal{T}_{h}(\Omega)$, then we
denote $\Vh(\omega):=\{ \varphi\vert_{\omega},\;\varphi\in
\Vh(\Omega)\}$ and also consider finite element spaces on boundaries
$\Vh(\partial\omega):=\{ \varphi\vert_{\partial\omega},\;\varphi\in
\Vh(\Omega)\}$.

\quad\\
We are interested in domain decomposition by substructuration, which
leads to introducing a non-overlapping subdomain partition of the
computational domain.
\begin{equation}\label{SubdomainPartition}
  \begin{array}{ll}
    \overline{\Omega} = \cup_{j=1}^{\mJ}\overline{\Omega}_{j},
    & \text{with}\quad\Omega_{j}\cap\Omega_{k} = \emptyset\quad
      \text{for}\;j\neq k\\[5pt]
    \Sigma := \cup_{j=1}^{\mJ}\Gamma_{j},
    & \text{where}\quad
        \Gamma_{j}:=\partial\Omega_{j},
  \end{array}
\end{equation}
where each $\Omega_{j}\subset \Omega$ is itself a polyhedral domain
that is exactly resolved by the triangulation. We do not make any
further assumption regarding the subdomain partitioning. In
accordance with the notations of the previous section, we introduce
\begin{equation}\label{SkeletonTraceMap}
  \begin{aligned}
    & \Vh(\Sigma):=\mrm{Im}(\mathscr{B}) =\{v\vert_{\Sigma}: v\in\Vh(\Omega)\}\\
    & \text{where}\;\;\mathscr{B}(u):= u\vert_{\Sigma}.
  \end{aligned}
\end{equation}
The space $\Vh(\Sigma)$ consists in (single valued) finite element
functions defined over the skeleton that is a surface with multiple
branches i.e. the union of all interfaces which is \textit{neither a
  boundary, nor even a manifold}.  By construction
$\mathscr{B}:\Vh(\Omega)\to \Vh(\Sigma)$ is surjective. Next we
introduce continuous and discrete function spaces naturally associated
to the multi-domain setting
\begin{equation}\label{ContinuousNorms}
  \begin{aligned}
    & \mbVh(\Omega) := \Vh(\Omega_1)\times \dots\times \Vh(\Omega_\mJ),\\
    & \mbVh(\Sigma) := \Vh(\Gamma_{1})\times\cdots\times \Vh(\Gamma_{\mJ}).
  \end{aligned}
\end{equation}
Since they are cartesian products, these spaces are made of tuples of
(volume based) functions.  The "broken space" $\mbVh(\Omega)$ is
naturally identified with those functions that are piecewise
$\mathbb{P}_1$-Lagrange in each subdomain. From this perspective, the
space $\Vh(\Omega)$ is embedded into $\mbVh(\Omega)$ by means of the
embedding operator $\mathscr{R}:\Vh(\Omega)\to \mbVh(\Omega)$ defined by
\begin{equation}\label{VolumeRestriction}
  \mathscr{R}(u):=(u\vert_{\Omega_1},\dots, u\vert_{\Omega_\mJ}).
\end{equation}
The range space $\mrm{Im}(\mathscr{R})$ can be identified with those
functions that are globally $\mathbb{P}_1$-Lagrange in the whole
computational domain, including through interfaces $\Gamma_j\cap
\Gamma_k$. The space $\mbVh(\Sigma)$ can be obtained by taking
interior traces of functions belonging to $\mbVh(\Omega)$
resp. $\mbX_h(\Omega)$. This motivates the introduction of a
multi-domain trace map $\Bh:\mbVh(\Omega)\to \mbVh(\Sigma)$ defined by
\begin{equation}\label{TraceOperator}
  \Bh(v_1,\dots,v_\mJ):= (v_1\vert_{\Gamma_{1}},\dots,v_\mJ\vert_{\Gamma_{\mJ}})
\end{equation}
for $v = (v_1,\dots,v_\mJ)\in\mbVh(\Omega)$.  This trace operator
\eqref{TraceOperator} surjectively maps $\mbVh(\Omega)$ onto
$\mbVh(\Sigma)$. We emphasize that the boundary trace map
\eqref{TraceOperator} is subdomain-wise block-diagonal. Since we are
in a finite dimensional context, and $\overline{\Bh(v)} =
\Bh(\overline{v})\;\forall v\in \mbVh(\Omega)$, according to
\cite[Thm. 4.7 \& 4.12]{zbMATH01022519} we have $\mrm{Range}(\Bh^*) =
\mrm{Ker}(\Bh)^\circ := \{\phi\in \mbVh(\Omega)^*:\; \langle
\phi,v\rangle = 0\;\forall v\in \mrm{Ker}(\Bh)\}$. A tuple of traces
can also be obtained by taking restrictions of single valued function
defined on the skeleton, which motivates the introduction of the
restriction operator $\mR:\Vh(\Sigma)\to\mbVh(\Sigma)$ defined by 
\begin{equation}\label{SurfaceRestriction}
  \mrm{R}(p):=(p\vert_{\Gamma_1},\dots,p\vert_{\Gamma_\mJ}).
\end{equation}
We have introduced two trace maps $\mathscr{B}$ and $\mB$ respectively
defined by \eqref{SkeletonTraceMap} and \eqref{TraceOperator}, and two
restriction maps $\mathscr{R}$ and $\mR$ respectively defined by
\eqref{VolumeRestriction} and \eqref{SurfaceRestriction}. These
four operators obviously satisfy the following identity
\begin{equation}\label{CommutationRelation}
  \mB\cdot\mathscr{R} = \mR\cdot\mathscr{B}. 
\end{equation}

\section{Harmonic liftings}\label{HarmonicLifting}

Consider two positive measurable functions
$\alpha,\beta:\Omega\to (0,+\infty)$ with $\alpha_-\leq \alpha(\bx)\leq \alpha_+,
\beta_-\leq \beta(\bx)\leq \beta_+\forall\bx\in \Omega$ for fixed
constants $\alpha_{\pm},\beta_{\pm}>0$. We first focus on a discrete
operator $\mathscr{L}:\Vh(\Omega)\to\Vh(\Omega)^*$  defined
by
\begin{equation}\label{DiscreteVF}
  \langle \mathscr{L}(u),v\rangle
  :=\int_{\Omega}\alpha\nabla u\cdot\nabla v + \beta\, u v\; d\bx
\end{equation}
According to the assumptions on $\beta,\alpha$ this operator induces a
scalar product over $\Vh(\Omega)$ with attached norm $\Vert
u\Vert_{\mathscr{L}}^2 = \langle \mathscr{L}(u),\overline{u}\rangle$. As
$\mathscr{B}:\Vh(\Omega)\to \Vh(\Sigma)$ is surjective, see
\eqref{SkeletonTraceMap}, following the abstract framework of the appendix,
we can consider its Moore-Penrose pseudo-inverse
$\mathscr{B}^{\dagger}_{\mathscr{L}}:\Vh(\Sigma)\to\Vh(\Omega)$ defined through
\eqref{DefPseudoInverse} with respect to $\Vert \cdot\Vert_{\mathscr{L}}$.

\quad\\
On the other hand, the operator $\mathscr{L}$
can be decomposed in accordance with the non-overlapping subdomain
partition \eqref{SubdomainPartition}, which gives rise to a bounded
block-diagonal operator $\mL:\mbVh(\Omega)\to \mbVh(\Omega)^*$ defined
by
\begin{equation}
  \begin{aligned}
    & \mL = \mrm{diag}(\mL_{\Omega_1},\dots,\mL_{\Omega_\mJ})\\
    & \langle \mL_{\Omega_j}u,v\rangle:=\int_{\Omega_j}\alpha\nabla u\cdot\nabla v + \beta u v \,d\bx.
  \end{aligned}
\end{equation}
It is then clear from these definitions and \eqref{DiscreteVF} that
$\mathscr{L} = \mathscr{R}^*\mL\mathscr{R}$. The operator $\mL$
induces a scalar product over $\mbVh(\Omega)$ and we shall denote
$\Vert u\Vert_{\mL}^{2} := \langle \mL u,\overline{u}\rangle$ the
associated norm.  Since $\mB:\mbVh(\Omega)\to \mbVh(\Sigma)$ is
surjective, in accordance with the abstract framework of the appendix,
 we can consider its Moore-Penrose
pseudo-inverse $\mB^{\dagger}_{\mL}:\mbVh(\Sigma)\to \mbVh(\Omega)$
defined through \eqref{DefPseudoInverse} with respect to $\Vert
\cdot\Vert_{\mL}$.

\quad\\
By construction we have $\mB\cdot\mB^{\dagger}_{\mL}= \Id$.  Moreover
$\mB^{\dagger}_{\mL}(\mrm{Im}(\mR))\subset\mrm{Im}(\mathscr{R})$.  The
map $\mB^{\dagger}_{\mL}$ is itself subdomain-wise block diagonal.
Given a (tuple of) traces $p\in \mbVh(\Sigma)$, the tuple of functions
given by $\mB^{\dagger}_{\mL}(p)$ are commonly referred to as harmonic
liftings local to subdomains. The operator
$\mB^{\dagger}_{\mL}\cdot\mB:\mbVh(\Omega)\to \mbVh(\Omega)$ is an
$\mL$-orthogonal projector with kernel
$\mrm{Ker}(\mB^{\dagger}_{\mL}\cdot\mB) = \mrm{Ker}(\mB)$.  The next
proposition generalizes \eqref{CommutationRelation} by providing a
relation between pseudo-inverses and restriction operators.

\begin{prop}\label{CommutationPseudoInverses}\quad\\
  $\mathscr{R}\cdot\mathscr{B}^{\dagger}_{\mathscr{L}} = \mB^{\dagger}_{\mL}\cdot\mR$
\end{prop}
\noindent \textbf{Proof:}

Take an arbitrary $u\in \Vh(\Sigma)$ and set $\varphi:=
\mathscr{R}\cdot\mathscr{B}^{\dagger}_{\mathscr{L}}(u)\in
\mbVh(\Omega)$ and $r:= \mR(u)\in\mbVh(\Sigma)$.  Applying
\eqref{CommutationRelation} yields $\mB(\varphi) =
\mB\cdot\mathscr{R}\cdot\mathscr{B}^{\dagger}_{\mathscr{L}}(u) =
\mR\cdot\mathscr{B}\cdot\mathscr{B}^{\dagger}_{\mathscr{L}}(u) =
\mR(u) = r$. As a consequence, according to the characterization
\eqref{DefPseudoInverse} of $\mB^{\dagger}_{\mL}(r) = \mB^{\dagger}_{\mL}\cdot\mR(u)$, we have
\begin{equation*}
  \Vert \mB^{\dagger}_{\mL}\cdot\mR(u)\Vert_{\mL}\leq \Vert \varphi\Vert_{\mL}.
\end{equation*}
On the other hand we have $\mB\cdot\mB^{\dagger}_{\mL}(r) = r =
\mR(u)\in \mrm{Im}(\mR)$ and since, for $w\in \mbVh(\Omega)$,
$\mB(w)\in \mrm{Im}(\mR) \Rightarrow w\in \mrm{Im}(\mathscr{R})$ we
conclude that there exists $v\in \Vh(\Omega)$ such that
$\mB^{\dagger}_{\mL}(r) = \mathscr{R}(v)$.  Applying
\eqref{CommutationRelation} yields $\mR\cdot\mathscr{B}(v) =
\mB\cdot\mathscr{R}(v) = \mB\cdot\mB^{\dagger}_{\mL}(r) = r =
\mR(u)\Rightarrow \mathscr{B}(v)-u\in \mrm{Ker}(\mR) = \{0\}$ hence
$\mathscr{B}(v)= u$. According to the characterization
\eqref{DefPseudoInverse} of
$\mathscr{B}^{\dagger}_{\mathscr{L}}(u)$ we deduce
\begin{equation}
  \begin{aligned}
    \Vert \varphi\Vert_{\mL}
    =  \Vert \mathscr{R}\cdot\mathscr{B}^{\dagger}_{\mathscr{L}}(u)\Vert_{\mL}
    & =    \Vert \mathscr{B}^{\dagger}_{\mathscr{L}}(u)\Vert_{\mathscr{L}}\\
    & \leq \Vert v\Vert_{\mathscr{L}} = \Vert \mathscr{R}(v)\Vert_{\mL} =
    \Vert \mB^{\dagger}_{\mL}\cdot\mR(u)\Vert_{\mL}. 
  \end{aligned}
\end{equation}
In conclusion we have established that $\Vert \varphi\Vert_{\mL} =
\Vert \mB^{\dagger}_{\mL}\cdot\mR(u)\Vert_{\mL}$ but the unique solvability
of the minimization problem \eqref{DefPseudoInverse} proves that $\varphi =
\mB^{\dagger}_{\mL}\cdot\mR(u)$ which rewrites
$\mathscr{R}\cdot\mathscr{B}^{\dagger}_{\mathscr{L}}(u) = \mB^{\dagger}_{\mL}\cdot\mR(u)$.
Since $u$ was arbitrarily chosen in $\Vh(\Sigma)$ this finishes the proof.
\hfill $\Box$

\quad\\
In the context of substructuring domain decomposition we shall be
particularly interested in the so-called Schur complement operator
$\mT_{\mL}:\mbVh(\Sigma)\to \mbVh(\Sigma)^*$ defined by
\begin{equation}\label{DtNMap}
  \mT_{\mL}^{-1}:=\mB\mL^{-1}\mB^{*}.
\end{equation}
This should be understood as a discrete counterpart of a
(subdomain-wise block-diagonal) Dirichlet-to-Neumann map associated to
$\mL$. As both $\mL$ and $\mB$ are block-diagonal, so is $\mT_{\mL} =
\mrm{diag}(\mT_{\mL,\Omega_1},\dots, \mT_{\mL,\Omega_\mJ})$ where each
$\mT_{\mL,\Omega_j}:\Vh(\Gamma_j)\to \Vh(\Gamma_j)^*$.

\begin{thm}\label{FormulaInverse}\quad\\
  $(\mR^*\mT_{\mL}\mR)^{-1} = \mathscr{B}\mathscr{L}^{-1}\mathscr{B}^*$.
\end{thm}
\noindent \textbf{Proof:}

We need to study
$(\mR^*\mT_{\mL}\mR)(\mathscr{B}^*\mathscr{L}^{-1}\mathscr{B})$ and
prove that this expression reduces to $\Id$. Applying Lemma
\ref{TransfoInverse} we have $\mT_{\mL}=(\mB^{\dagger}_{\mL})^{*}\mL
\mB^{\dagger}_{\mL}$. Moreover we have $\mL (\mB^{\dagger}_{\mL} \mB)
= (\mB^{\dagger}_{\mL} \mB)^*\mL$ according to
\eqref{SelfAdjointProjector}. Combining these identities with
\eqref{CommutationRelation} leads to
\begin{equation}\label{CalculIntermediaire1}
  \begin{aligned}
    (\mR^*\mT_{\mL}\mR)(\mathscr{B}^*\mathscr{L}^{-1}\mathscr{B})
    & =  \mR^*(\mB^{\dagger}_{\mL})^{*}\mL \mB^{\dagger}_{\mL} \mR\mathscr{B}\mathscr{L}^{-1}\mathscr{B}^* \\
    & = \mR^*(\mB^{\dagger}_{\mL})^{*}\mL \mB^{\dagger}_{\mL} \mB\mathscr{R}\mathscr{L}^{-1}\mathscr{B}^*\\
    & = \mR^*(\mB^{\dagger}_{\mL})^{*}(\mB^{\dagger}_{\mL} \mB)^{*}\mL \mathscr{R}\mathscr{L}^{-1}\mathscr{B}^*
  \end{aligned}
\end{equation}
Next we have $\mB \mB^{\dagger}_{\mL} = \Id$ hence
$(\mB^{\dagger}_{\mL})^{*}(\mB^{\dagger}_{\mL}\mB)^{*} =
(\mB^{\dagger}_{\mL})^{*}$ which simplifies two factors in the
expression above. Taking the transpose of the identity in Proposition
\ref{CommutationPseudoInverses} yields $\mR^*(\mB^{\dagger}_{\mL})^{*}
= (\mathscr{B}^{\dagger}_{\mathscr{L}})^*\mathscr{R}^*$. As a
consequence, coming back to \eqref{CalculIntermediaire1} and taking
into account that $\mathscr{R}^*\mL \mathscr{R} = \mathscr{L}$, we
finally obtain
\begin{equation*}
  \begin{aligned}
    (\mR^*\mT_{\mL}\mR)(\mathscr{B}^*\mathscr{L}^{-1}\mathscr{B})
    & = \mR^*(\mB^{\dagger}_{\mL})^{*}\mL \mathscr{R}\mathscr{L}^{-1}\mathscr{B}^*\\
    & = (\mathscr{B}^{\dagger}_{\mathscr{L}})^*\mathscr{R}^*\mL \mathscr{R}\mathscr{L}^{-1}\mathscr{B}^*\\
    & = (\mathscr{B}^{\dagger}_{\mathscr{L}})^*\mathscr{L}\mathscr{L}^{-1}\mathscr{B}^* = 
    (\mathscr{B}^{\dagger}_{\mathscr{L}})^*\mathscr{B}^* = \Id.
  \end{aligned}
\end{equation*}
\hfill $\Box$

\quad\\
As direct byproduct, the preceding formula shows how any
preconditioner $\tilde{\mathscr{L}}^{-1}$ for $\mathscr{L}$ gives rise
to a preconditioner for the Schur complement system
$\mR^*\mT_{\mL}\mR$.

\begin{cor}\label{PrecondtitioningProcess}\quad\\
  For any hermitian positive definite linear map
  $\tilde{\mathscr{L}}:\Vh(\Omega)\to \Vh(\Omega)^{*}$, we have
  \begin{equation*}
    \mrm{cond}( (\mathscr{B}\tilde{\mathscr{L}}^{-1}\mathscr{B}^*) (\mR^*\mT_{\mL}\mR) )\leq
    \mrm{cond}(\tilde{\mathscr{L}}^{-1}\mathscr{L} )
  \end{equation*}
\end{cor}
\noindent \textbf{Proof:}

According to the formula of Theorem \ref{FormulaInverse}, we have  the identity 
$\mrm{cond}( (\mathscr{B}\tilde{\mathscr{L}}^{-1}\mathscr{B}^*) (\mR^*\mT_{\mL}\mR) ) = \lambda_+/\lambda_{-}$
where the upper extremal eigenvalue $\lambda_+$ is defined by
\begin{equation*}
\lambda_{+} :=
\sup_{u\in \Vh(\Sigma)^{*}\setminus\{0\}}
\frac{\langle (\mathscr{B}\tilde{\mathscr{L}}^{-1}\mathscr{B}^*) u,\overline{u}\rangle}{
  \langle(\mathscr{B}\mathscr{L}^{-1}\mathscr{B}^*) u,\overline{u}\rangle}
\end{equation*}
and the lower eigenvalue $\lambda_-$ is defined similarly taking the infimum
instead of the supremum. Because $\mathscr{B}^*(\Vh(\Sigma)^*)\subset \Vh(\Omega)^*$, the
upper eigenvalue is readily bounded by
\begin{equation*}
  \lambda_{+}\leq \sup_{u\in \Vh(\Omega)^{*}\setminus\{0\}}
  \frac{\langle \tilde{\mathscr{L}}^{-1}(u),\overline{u}\rangle}{
    \langle\mathscr{L}^{-1}(u),\overline{u}\rangle} =
  \sup_{v\in \Vh(\Omega)\setminus\{0\}}
  \frac{\langle \mathscr{L}(\tilde{\mathscr{L}}^{-1}\mathscr{L})v,\overline{v}\rangle}{
  \langle\mathscr{L}v,\overline{v}\rangle} = \sup \mathfrak{S}(\tilde{\mathscr{L}}^{-1}\mathscr{L})
\end{equation*}
where $\mathfrak{S}(\tilde{\mathscr{L}}^{-1}\mathscr{L})$ refers to the spectrum of
$\tilde{\mathscr{L}}^{-1}\mathscr{L}$. In the same manner, we have
$\lambda_-\geq \inf \mathfrak{S}(\tilde{\mathscr{L}}^{-1}\mathscr{L})$
so that $\lambda_+/\lambda_-\leq \sup \mathfrak{S}(\tilde{\mathscr{L}}^{-1}\mathscr{L})/
\inf \mathfrak{S}(\tilde{\mathscr{L}}^{-1}\mathscr{L}) = \mrm{cond}(\tilde{\mathscr{L}}^{-1}\mathscr{L})$.
\hfill $\Box$

\section{Neumann-Neumann preconditioner}\label{RappelNeumannNeumann}

We shall now describe in detail the so-called Neumann-Neumann method,
a well-established preconditioning strategy for the scalar Schur
complement system.  We briefly recall the principle of this
substructuring technique that we recast in the notations we have
introduced thus far. This preconditioner is discussed at length
e.g. in \cite[\S 4.3.2]{MR1410757}, \cite[\S 3.7.1]{MR2445659},
\cite[Chap.6]{MR2104179} or \cite[\S 2.3]{MR3013465}.

\quad\\
First of all observe that $\mR:\Vh(\Sigma)\to\mbVh(\Sigma)$ defined by
\eqref{SurfaceRestriction} is injective by construction, so we can apply
Lemma \ref{TransfoInverse2} to characterize the inverse of the Schur
complement system $\mR^*\mT_{\mL}\mR$ which yields
\begin{equation}\label{IdealNNPrec}
    (\mR^*\mT_{\mL}\mR)^{-1} = (\mR^{\dagger}_{\mT_{\mL}})\mT_{\mL}^{-1}(\mR^{\dagger}_{\mT_{\mL}})^{*}\\
\end{equation}
where the right hand side involves the pseudo-inverse of $\mR$ with
respect to the scalar product induced by $\mT_{\mL}$.  The idea of the
Neumann-Neumann method consists in considering a pseudo-inverse for
$\mR$ with respect to a different, more convenient scalar product i.e.
we approximate $(\mR^*\mT_{\mL}\mR)^{-1}$ with
$\mrm{P}_{\textsc{nn}}:=\mR^{\dagger}_{\mD}\mT_{\mL}^{-1}(\mR^{\dagger}_{\mD})^{*}$ and
$\mR^{\dagger}_{\mD}:\mbVh(\Sigma)\to \Vh(\Sigma)$ refers to the
pseudo-inverse of $\mR$ with respect to the scalar product induced by
$\mD:\mbVh(\Sigma)\to \mbVh(\Sigma)^*$ defined by
\begin{equation*}
  \begin{aligned}
    & \mD = \mrm{diag}(\mD_{\Gamma_1},\dots,\mD_{\Gamma_\mJ})\quad \text{with}\\
    & \langle \mD_{\Gamma_j}(u),v\rangle := \sum_{\bx\in \mrm{vtx}(\Gamma_j)} u(\bx)v(\bx).
  \end{aligned}
\end{equation*}
where $\mrm{vtx}(\Gamma_j) = \{$vertices of the mesh located on
$\Gamma_j \}$.  Each local operator $\mD_{\Gamma_j}$ is naturally
represented as the identity matrix in the standard nodal basis of
$\mathbb{P}_1$-shape functions of $\Vh(\Gamma_j)$.  This is why
$\mD_{\Sigma}:=\mR^*\mD\mR$ is trivially invertible, being itself
represented by a diagonal matrix in the nodal basis of
$\mathbb{P}_1$-shape functions of $\Vh(\Sigma)$ according to the
expression
\begin{equation*}
  \begin{aligned}
    & \langle \mD_{\Sigma}(u),v\rangle = \sum_{\bx\in \mrm{vtx}(\Sigma)} \mrm{deg}(\bx) u(\bx)v(\bx)\\
    & \text{where}\;\;\mrm{deg}(\bx) := \mrm{card}\{j\in \{1,\dots,\mJ\},\bx\in\Gamma_j\}.
  \end{aligned}
\end{equation*}
The pseudo-inverse of the operator $\mR$ relative to the scalar
product $\mD$ is then given by the expression $\mR_{\mD}^{\dagger} =
\mD_{\Sigma}^{-1}\mR^{*}\mD$ and the Neumann-Neumann preconditioner
$\mrm{Q}_{\textsc{nn}}:\Vh(\Sigma)^*\to\Vh(\Sigma)$ for the Schur
complement system $\mR^*\mT_\mL\mR$ is then given by
\begin{equation}\label{NeumannNeumannPrec}
  \begin{aligned}
    \mrm{Q}_{\textsc{nn}}
    & := (\mR_{\mD}^{\dagger})\mT_{\mL}^{-1}(\mR_{\mD}^{\dagger})^{*}\\
    & \phantom{:}=\mD_{\Sigma}^{-1}\mR^{*}\mD\mT_{\mL}^{-1}\mD\mR\mD_{\Sigma}^{-1}.
  \end{aligned}
\end{equation}
We see that \eqref{NeumannNeumannPrec} only differs from \eqref{IdealNNPrec}
by $\mR_{\mD}^{\dagger}$ replacing $\mR_{\mT_{\mL}}^{\dagger}$. The performance analysis
of this preconditioner is rather standard. For a fixed subdomain partition,
and under technical assumptions, the condition number is proved to only deteriorate
logarithmically with respect to meshwidth, see e.g. \cite[\S 2.6]{MR3013465}.

\section{Numerical illustration}\label{NumericIllustrationSec}

In this section, we shall present numerical results confirming the
relevancy of \eqref{NeumannNeumannPrec} as a preconditioner. We
consider the computational domain $\Omega = (0,1)\times(0,1)\times
(0,1)\subset \mathbb{R}^{3}$ regularly subdivided in $3\times 3\times
3 = 27$ subdomains. We generate a mesh of $\Omega$ conforming with
this subdivision, see Figure \ref{Fig1} below. Mesh generation is
achieved by means of \texttt{gmsh} \cite{MR2566786}.
\begin{figure}[h!]
\centering\includegraphics[height=5cm,trim=1000 50 1000 50]{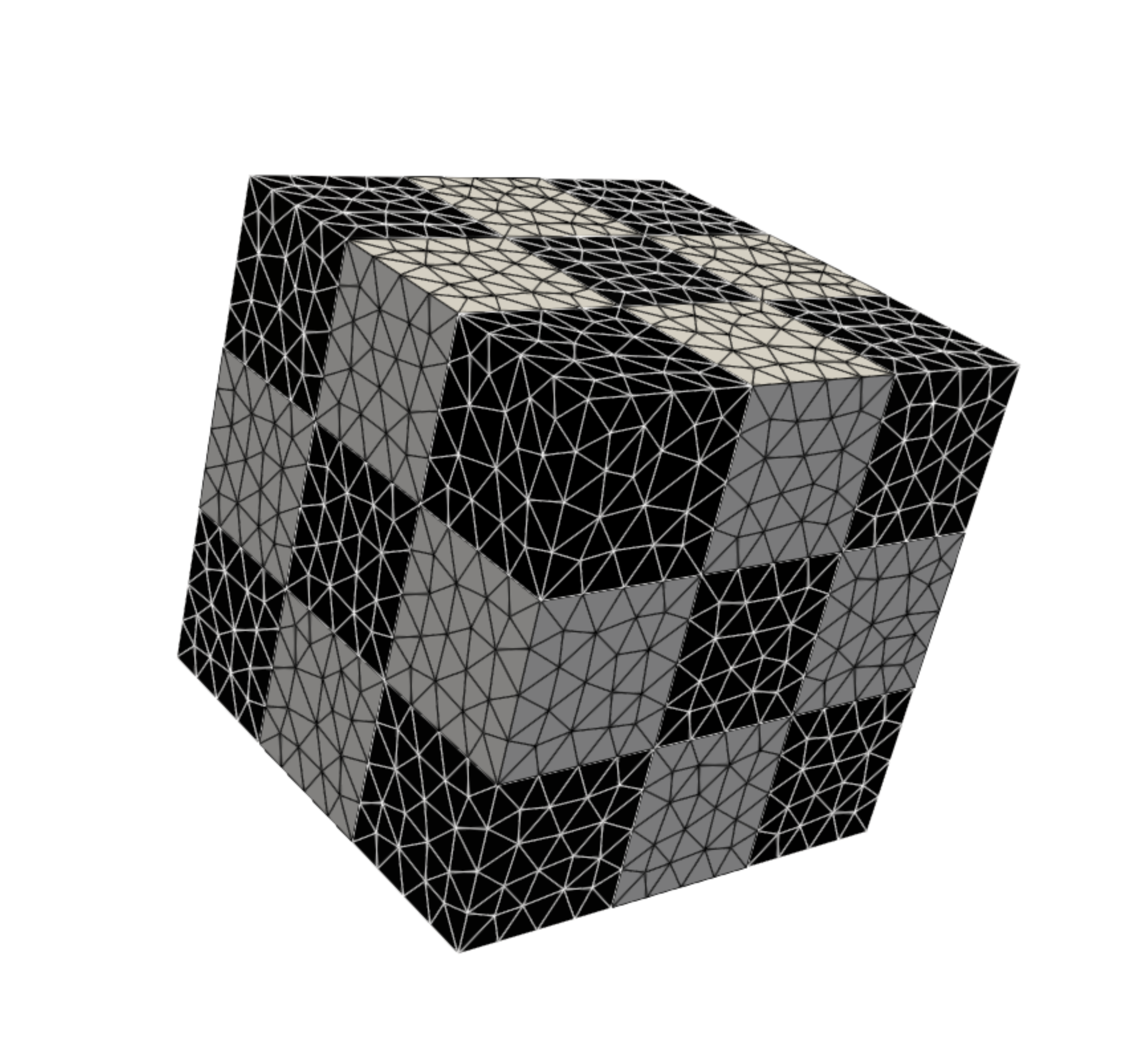}
\caption{Global computational domain}\label{Fig1}
\end{figure}
We take $\alpha = \beta = 1$. First we examine the performance of the
scalar Neumann-Neumann preconditioner \eqref{NeumannNeumannPrec} for
this actual setup: a vector $\bu_{\textsc{ex}} \in \Vh(\Sigma)$ is
drawn randomly, set $\bff = (\mR^*\mT_\mL\mR)\bu_{\textsc{ex}}$ and
then solve the linear system
$\mrm{Q}_{\textsc{nn}}(\mR^*\mT_\mL\mR)\bu =
\mrm{Q}_{\textsc{nn}}(\bff)$ by means of a Preconditioned Conjugate
Gradient (PCG) algorithm \cite[\S 9.2]{zbMATH01953446}.  We conduct
this numerical solve for 5 different meshes. In the figure below we
plot the convergence history and record the required number of
iteration to reach $10^{-9}$ relative residual for each mesh.

\noindent 
\begin{minipage}{0.45\linewidth}
    \centering\includegraphics[height=6cm]{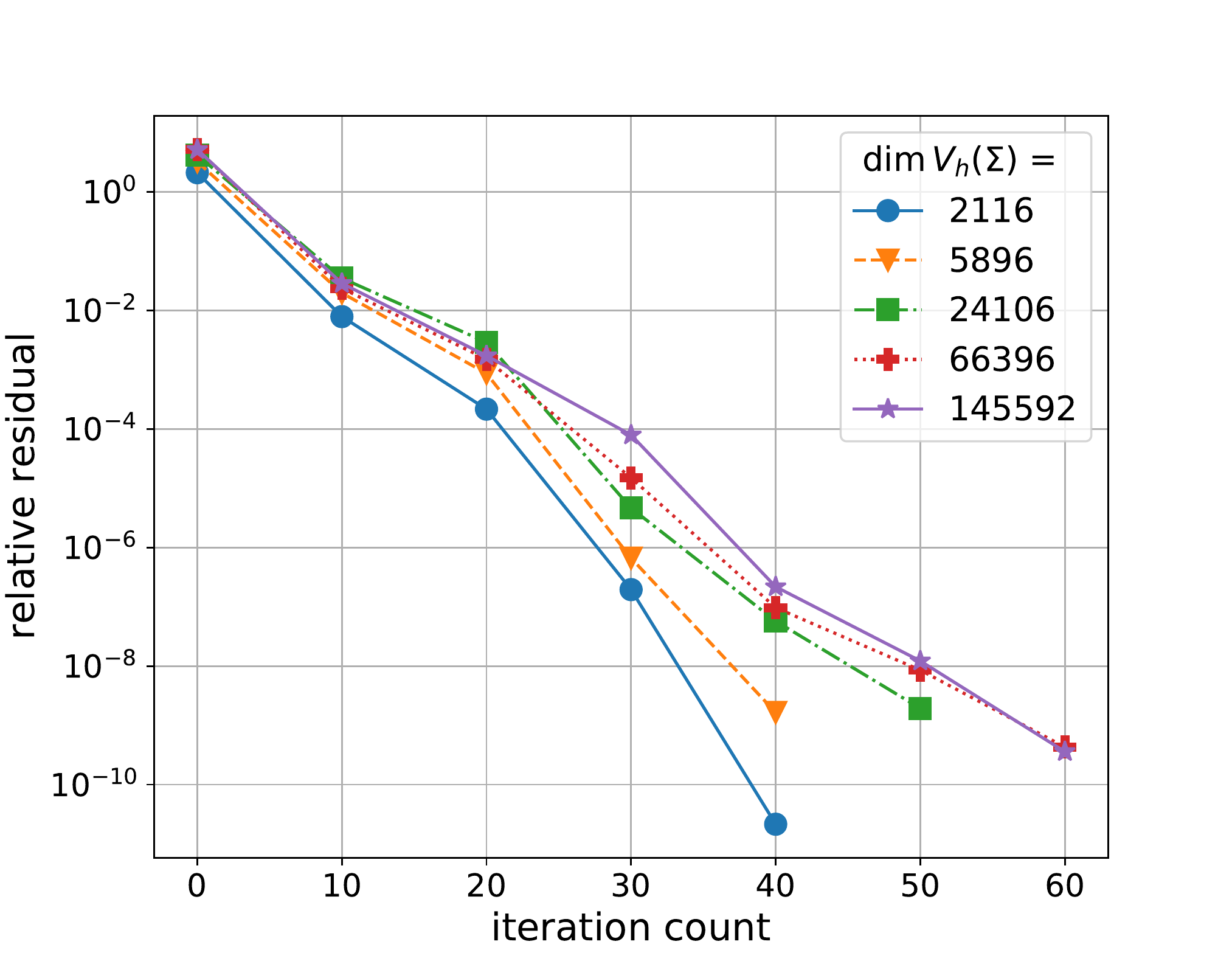}
\end{minipage}
\hspace{2cm}
\begin{minipage}{0.4\linewidth}
  \centering 
  \begin{tabular}{r|r|p{1cm}}
    $\mrm{dim}\,\Vh(\Sigma)$ & $\mrm{dim}\,\Vh(\Omega)$ & \# iter\\
    \hline
    2116   \hspace{0.1cm} & 2364   \hspace{0.1cm} &\hspace{0.3cm} 41\\
    5896   \hspace{0.1cm} & 8700   \hspace{0.1cm} &\hspace{0.3cm} 48\\
    24106  \hspace{0.1cm} & 61093  \hspace{0.1cm} &\hspace{0.3cm} 56\\
    66396  \hspace{0.1cm} & 258046 \hspace{0.1cm} &\hspace{0.3cm} 62\\
    145592 \hspace{0.1cm} & 805896 \hspace{0.1cm} &\hspace{0.3cm} 63
  \end{tabular}\quad\\[5pt]
  \captionof{table}{Required iterations to reach $10^{-9}$ relative residual
  in PCG applied to $\mR^*\mT_\mL\mR$ preconditioned
  with $\mrm{Q}_{\textsc{nn}}$.} \label{tab:title} 
  
\end{minipage}

\quad\\[10pt]
Here $\mrm{dim}(\Vh(\Sigma))$ is the number of unknowns
i.e.  the size of the linear system. The growth of the required number
of iterations is mild compared to the growth of the number of unknowns.

\section{Edge element discretization}
The results of Section \ref{HarmonicLifting} provide guidelines for converting a
preconditioner for a volume based PDE into a preconditioner for the
corresponding Schur complement system. This is the message underlying
Corollary \ref{PrecondtitioningProcess}. We wish to examine how this applies to the
Hiptmair-Xu preconditioner.

\quad\\
We introduce vector counterparts of many of the definitions introduced
in Section \ref{GeometricConfig} and, for some of these, adopt the same notation
which complies with consistency. First we consider volume based edge
finite element space
\begin{equation}\label{EdgeElt}
  \begin{aligned}
    \Wh(\Omega) := \{\bu\in \mathbf{H}(\boldsymbol{curl},\Omega),\;
    \forall \tau\in \mathcal{T}_{h}(\Omega)\;\exists \boldsymbol{a}_\tau,\boldsymbol{b}_\tau\in \mathbb{C}^3\;\\
    \bu(\bx) = \boldsymbol{a}_\tau + \boldsymbol{b}_\tau\times \bx,\;\forall \bx\in \tau\}
  \end{aligned}
\end{equation}
These are lowest order curl-conforming edge elements, see \cite[\S 5.5]{MR2059447}.
Denoting $\mrm{edge}(\Omega)$ the collection of edges of the triangulation
$\mathcal{T}_h(\Omega)$, and letting $\tau_e\in \RR^3$ refer to a unit vector
tangent to $e\in \mrm{edge}(\Omega)$, the dual edge element space is
spanned by $\Wh(\Omega)^* = \mrm{span}\{\phi_e^*\}_{e\in \mrm{edge}(\Omega)}$
with degrees of freedom $\phi_e^*:\mathscr{C}^{0}(\overline{\Omega})^{3}\to \mathbb{C}$ defined by
\begin{equation*}
  \langle \phi_{e}^*,\bv\rangle:=\int_e\bv\cdot \tau_e d\sigma\quad \forall \bv\in \mathscr{C}^{0}(\overline{\Omega})^3.
\end{equation*}
The shape function $\phi_e$ are defined as the unique element of $\Wh(\Omega)$
satisfying $\langle \phi_{e}^*,\phi_e\rangle = 1$ and $\langle \phi_{e}^*,\phi_{e'}\rangle = 0$
for $e,e'\in \mrm{edge}(\Omega)$, $e\neq e'$.
Associated to this edge element space, one classically define a local
N\'ed\'elec interpolation operator $\Pi_{\Omega}:\mathscr{C}^{0}(\overline{\Omega})^3\to
\Wh(\Omega)$ characterized by the property that for all $\bu\in
\mathscr{C}^{0}(\overline{\Omega})^3$, the vector field
$\Pi_{\Omega}(\bu)\in \Wh(\Omega)$ is the unique edge element vector
field satisfying
\begin{equation*}
    \Pi_{\Omega}(\bu) = \sum_{e\in\mrm{edge}(\Omega)}\langle \phi_e^*,\bu\rangle \phi_e.
\end{equation*}
Setting $\Wh(\Omega_j):=\{ \bu\vert_{\Omega_j},\bu\in \Wh(\Omega)\}$
and $\mbWh(\Omega):= \Wh(\Omega_1)\times\dots\times \Wh(\Omega_\mJ)$, 
the definition of the restriction operator $\underline{\mathscr{R}}:\Wh(\Omega)\to
\mbWh(\Omega)$ naturally extends to the vector valued context setting
$\underline{\mathscr{R}}(\bu) := (\bu\vert_{\Omega_1},\dots,\bu\vert_{\Omega_\mJ})$.

\quad\\
We also need to consider surface vector fields and edge elements
tangent to the skeleton of the subdomain partition. Let $\bn_{j}$ be
the unit normal vector field at $\Gamma_j$ directed toward the
exterior of $\Omega_j$, and define $\bn_{\Sigma}:=\bn_j$ on
$\Gamma_j\cap \Gamma_k, j<k$ and $\bn_{\Sigma}:=\bn_j$ on
$\Gamma_j\cap\partial\Omega$. Then we set
\begin{equation}\label{SurfaceEdgeElt}
  \begin{aligned}
    \mathscr{C}^{0}_{\textsc{t}}(\Sigma)
    & :=\{\bn_{\Sigma}\times
    \bu\vert_{\Sigma}\times\bn_{\Sigma},\;\bu\in \mathscr{C}^0(\Omega)^{3}\}\\
    \Wh(\Sigma)
    & := \{\bn_{\Sigma}\times \bu\vert_{\Sigma}\times\bn_{\Sigma},\;\bu\in \Wh(\Omega)\;\}
  \end{aligned}
\end{equation}
where $\bn_{\Sigma}\times \bu\vert_{\Sigma}\times\bn_{\Sigma}$ is
simply the tangent part of $ \bu\vert_{\Sigma}$. The space
$\Wh(\Sigma)$ is properly defined because in $\Wh(\Omega)$, by
construction, tangent traces match across any mesh face common to two
neighboring elements. The definition of the skeleton trace operator
is then naturally extended by
\begin{equation*}
  \underline{\mathscr{B}}(\bu) :=  \bn_{\Sigma}\times \bu\vert_{\Sigma}\times\bn_{\Sigma}
\end{equation*}
so that, by construction, $\underline{\mathscr{B}}$ surjectively maps
$\Wh(\Omega)$ (resp.~$\mathscr{C}^0(\Omega)^{3}$) onto $\Wh(\Sigma)$
(resp.~$\mathscr{C}^{0}_{\textsc{t}}(\Sigma)$).  This operator does
not depend on the actual direction of the normal vector field as
$\bn_{\Sigma}$ comes into play twice in the formula above.  Setting
$\Wh(\Gamma_j):=\{\bn_j\times
\bu\vert_{\Gamma_j}\times\bn_j,\;\bu\in\Wh(\Omega_j)\}$ and
$\mbWh(\Sigma) = \Wh(\Gamma_1)\times \dots\times \Wh(\Gamma_\mJ)$, the
definition of the local surface restriction operator $\underline{\mR}:\Wh(\Sigma)
\to \mbWh(\Sigma)$ is naturally extended by
$\underline{\mR}(\bu):=(\bu\vert_{\Gamma_1},\dots,\bu\vert_{\Gamma_\mJ})$ for all
$\bu\in \Wh(\Sigma)$. Finally we also define a subdomain-wise boundary
trace operator
\begin{equation}
  \underline{\mB}(\bu_1,\dots,\bu_\mJ):=(\bn_j\times\bu_j\vert_{\Gamma_j}\times \bn_j)_{j=1,\dots,\mJ}
\end{equation}
for $\bu_j\in\Wh(\Omega_j), j=1\dots \mJ$, which is a natural vector
extension of \eqref{TraceOperator}. The operator above induces a
surjection $\underline{\mB}:\mbWh(\Omega)\to\mbWh(\Sigma)$. Again, the
vector trace maps satisfy a commutativity relation with restriction
operators $\underline{\mB}\cdot\underline{\mathscr{R}} =
\underline{\mR}\cdot\underline{\mathscr{B}}$.
The gradient operators induce operators
both in the volume $\mrm{G}_{\Omega}:\Vh(\Omega)\to \Wh(\Omega)$ and on the
skeleton $\mrm{G}_{\Sigma}:\Vh(\Sigma)\to \Wh(\Sigma)$, and defined by
\begin{equation*}
  \mrm{G}_{\Omega}(u) := \nabla u,\quad \mrm{G}_{\Sigma}(v) := \nabla_{\Sigma} v
\end{equation*}
for $u\in \Vh(\Omega)$ and any $v\in\Vh(\Sigma)$, where
$\nabla_{\Sigma}$ refers to a surface gradient tangent to $\Sigma$. These also
satisfy commutation relation with the boundary trace operator
$\mrm{G}_{\Omega}\cdot \mathscr{B} = \underline{\mathscr{B}}\cdot \mrm{G}_{\Sigma}$.

\quad\\
Relying on the Lagrange interpolation operator
$\Pi_\Omega:\mathscr{C}^0(\overline{\Omega})^{3}\to \Wh(\Omega)$, any
vector $\be\in \mathbb{C}^3$ induces a linear operator
$\Pi_\Omega^{\be}:\mathscr{C}^0(\overline{\Omega})\to\Wh(\Omega)$
defined the following formula
\begin{equation}
  \Pi_\Omega^{\be}(u):=\Pi_\Omega(\be\,u).
\end{equation}
Observe that, for an edge included in the skeleton $e\subset \Sigma$,
there exists a functional
$\varphi_e^{*}:\mathscr{C}^0_{\textsc{t}}(\overline{\Omega})\to
\mathbb{C}$ such that $\langle \varphi_e^{*},\underline{\mathscr{B}}(\bv)\rangle =
\langle \phi_e^{*},\bv\rangle$. For an edge on the skeleton, we can
set $\varphi_e:=\underline{\mathscr{B}}(\phi_e)$ so that $\mW_h(\Sigma) =
\mrm{span}_{e\in\mathcal{E}_h(\Sigma)}\{\varphi_e\}$, and we can define
a surface interpolation operator $\Pi_{\Sigma}:
\mathscr{C}^{0}_{\textsc{t}}(\Sigma)\to \Wh(\Sigma)$ by the formula
\begin{equation*}
  \begin{aligned}
    & \underline{\mathscr{B}}\cdot\Pi_{\Omega} = \Pi_{\Sigma}\cdot\underline{\mathscr{B}}\quad \text{where}\\
    & \Pi_{\Sigma}(\bu):=\sum_{e\in\mrm{edge}(\Sigma)}\langle \varphi_e^*,\bu\rangle \varphi_e
  \end{aligned}
\end{equation*}
We can also define a surface counterpart of $\Pi_{\Omega}^{\be}$ in a
similar manner. For an arbitrary
$u\in\mathscr{C}^{0}(\overline{\Omega})$, observe that
$\underline{\mathscr{B}}\cdot\Pi_{\Omega}^{\be}(u)
=\underline{\mathscr{B}}\cdot\Pi_\Omega(\be\,u) =
\Pi_\Sigma(\underline{\mathscr{B}}(\be)\,\mathscr{B}(u))$.
This leads to defining $\Pi_{\Sigma}^{\be}:\mathscr{C}^{0}(\Sigma)
\to \mathbb{C}$, for any $u\in\mathscr{C}^{0}(\Sigma)$, by the
expression
\begin{equation}
  \Pi_{\Sigma}^{\be}(u):=\Pi_\Sigma(\underline{\mathscr{B}}(\be)\,u).
\end{equation}

\section{Hiptmair-Xu preconditioner}\label{HXPrecSec}
Now for some constant coefficient $\gamma>0$, consider the
operator $\mathscr{M}:\Wh(\Omega)\to\Wh(\Omega)^*$ defined,
for any $\bu,\bv\in\Wh(\Omega)$ by the formula
\begin{equation}\label{GlobalMaxwellOperator}
  \langle \mathscr{M}(\bu),\bv\rangle:=\int_{\Omega}
  \boldsymbol{curl}(\bu)\cdot\boldsymbol{curl}(\bv) + \gamma^2 \bu\cdot\bv d\bx
\end{equation}
This operator induces a scalar product over $\Wh(\Omega)$ whose norm
is obviously equivalent to the standard norm over
$\mathbf{H}(\boldsymbol{curl},\Omega)$. We need to consider a
so-called Jacobi smoother, that consists in the diagonal part of this
operator and defined by
\begin{equation}
  \langle \mrm{diag}(\mathscr{M})\bu,\bv\rangle =
  \sum_{e\in\mathcal{E}_h(\Omega)} \langle \mathscr{M}(\phi_e),
  \phi_e\rangle \phi_e^*(\bu)\phi_e^*(\bv)
\end{equation}
Preconditioning $\mathscr{M}$ is challenging due to the term
$\boldsymbol{curl}\cdot \boldsymbol{curl}$ admitting a large kernel.
Hiptmair and Xu \cite{MR2361899} (see in particular Formula (7.3) in
this article) proposed a nodal preconditioner
$\tilde{\mathscr{M}}^{-1}:\Wh(\Omega)^*\to \Wh(\Omega)$ based on
regular decompositions of the fields \cite{MR4143285}. It takes the
following form
\begin{equation}\label{HXPrec}
  \tilde{\mathscr{M}}^{-1} = \mrm{diag}(\mathscr{M})^{-1}
    +(\mrm{G}_{\Omega})\mathscr{L}^{-1}(\mrm{G}_{\Omega})^*
    +\sum_{j=1,2,3}(\Pi_\Omega^{\be_j})\mathscr{L}^{-1}(\Pi_\Omega^{\be_j})^*.
\end{equation}
In this expression $\be_j,j=1,2,3$ refer to the vectors of the
cartesian canonical basis of $\RR^3$. In practice $\mathscr{L}^{-1}$ is
approximated by a preconditioner, see
\cite[Cor.2.3]{MR2361899}. Dealing with \eqref{GlobalMaxwellOperator}
by means of substructuration leads to considering the Schur complement
system $\underline{\mR}^*\mT_{\mrm{M}}\underline{\mR}: \Wh(\Sigma)\to
\Wh(\Sigma)^*$ where $\mrm{M}:\mbWh(\Omega)\to \mbWh(\Omega)^*$ is the
subdomain-wise block diagonal counterpart of the global operator
$\mathscr{M}$ i.e.
\begin{equation*}
  \begin{aligned}
    & \mrm{M}:=\mrm{diag}(\mrm{M}_{\Omega_1},\dots,\mrm{M}_{\Omega_\mJ})\\
    & \langle \mrm{M}_{\Omega_j}(\bu),\bv\rangle:=\int_{\Omega_j}
      \boldsymbol{curl}(\bu)\cdot\boldsymbol{curl}(\bv) + \gamma^2 \bu\cdot\bv\;d\bx.
  \end{aligned}
\end{equation*}
Moreover a counterpart of Theorem \ref{FormulaInverse} can be also
established in the case of the operator $\mathscr{M}$, following the
very same proof pattern, which yields the identity
$(\underline{\mR}^*\mT_{\mrm{M}}\underline{\mR})^{-1} =
\underline{\mathscr{B}}\mathscr{M}^{-1}\underline{\mathscr{B}}^{*}$.
Combining this formula with \eqref{HXPrec}, and taking advantage of
the commutation properties satisfied by the boundary trace and the
gradient operators, leads to a preconditioner for the Schur complement
system $\underline{\mR}^*\mT_{\mrm{M}}\underline{\mR}$
\begin{equation}\label{RestructuredMaxwellPrec}
  \begin{aligned}
  & \;\underline{\mathscr{B}}\tilde{\mathscr{M}}^{-1}\underline{\mathscr{B}}^{*}\\[10pt]
  & = \underline{\mathscr{B}}\,\mrm{diag}(\mathscr{M})^{-1}\underline{\mathscr{B}}^{*}
  +\underline{\mathscr{B}}\,\mrm{G}_{\Omega}\mathscr{L}^{-1}\mrm{G}_{\Omega}^*\,\underline{\mathscr{B}}^*
  +\sum_{j=1,2,3}\underline{\mathscr{B}}(\Pi_\Omega^{\be_j})
  \mathscr{L}^{-1}(\Pi_\Omega^{\be_j})^*\underline{\mathscr{B}}^*\\
  & = \underline{\mathscr{B}}\,\mrm{diag}(\mathscr{M})^{-1}\underline{\mathscr{B}}^{*}
  +\mrm{G}_{\Sigma}(\mathscr{B}\mathscr{L}^{-1}\mathscr{B}^*)\mrm{G}_{\Sigma}^*
  +\sum_{j=1,2,3}(\Pi_\Sigma^{\be_j})(\mathscr{B}\mathscr{L}^{-1}\mathscr{B}^*)(\Pi_\Sigma^{\be_j})^*\\
  & = \underline{\mathscr{B}}\,\mrm{diag}(\mathscr{M})^{-1}\underline{\mathscr{B}}^{*}
  +\mrm{G}_{\Sigma}(\mR^*\mT_{\mL}\mR)^{-1}\mrm{G}_{\Sigma}^*
  +\sum_{j=1,2,3}(\Pi_\Sigma^{\be_j})(\mR^*\mT_{\mL}\mR)^{-1}(\Pi_\Sigma^{\be_j})^*\\
  \end{aligned}
\end{equation}
Let us plug \eqref{NeumannNeumannPrec} into \eqref{RestructuredMaxwellPrec},
replacing $(\mR^*\mT_{\mL}\mR)^{-1}$ with its Neumann-Neumann preconditioner
$\mrm{Q}_{\textsc{nn}}$. We obtain the following expression of a substructured variant of
the Hiptmair-Xu preconditioner
\begin{equation}\label{FinalPrec}
  \mrm{Q}_{\textsc{hx}}:=\underline{\mathscr{B}}\mrm{diag}(\mathscr{M})^{-1}
  \underline{\mathscr{B}}^{*} +\mrm{G}_{\Sigma}\,\mrm{Q}_{\textsc{nn}}\,\mrm{G}_{\Sigma}^*
  +\sum_{j=1,2,3}(\Pi_\Sigma^{\be_j})\mrm{Q}_{\textsc{nn}}(\Pi_\Sigma^{\be_j})^*
\end{equation}
The next lemma shows that the performance of this preconditioner can
be assessed by combining the condition number estimates available for
the HX-preconditioner $\tilde{\mathscr{M}}^{-1}$ on the one hand, see
\cite[Thm.7.1]{MR2361899}, and for the Neumann-Neumann preconditioner
on the other hand, see \cite[\S 2.6]{MR3013465}. 

\begin{lem}\label{FinalCondEstimate}
  \begin{equation*}
    \mrm{cond}(\mrm{Q}_{\textsc{hx}}\underline{\mR}^*\mT_{\mrm{M}}\underline{\mR})
    \leq \mrm{cond}(\mrm{Q}_{\textsc{nn}} \mR^*\mT_{\mL}\mR)\cdot\mrm{cond}(\tilde{\mathscr{M}}^{-1}\mathscr{M}).
  \end{equation*}    
\end{lem}
\noindent \textbf{Proof:}

The proof is a direct consequence of the inequality
$\lambda_{-}\langle \mrm{Q}_{\textsc{nn}}(\bu), \overline{\bu}\rangle \leq
\langle (\mR^*\mT_{\mL}\mR)^{-1}\bu, \overline{\bu}\rangle\leq
\lambda_{+}\langle \mrm{Q}_{\textsc{nn}}(\bu), \overline{\bu}\rangle$
for all $\bu\in \Vh(\Sigma)^{*}$, where $\lambda_{\pm}$ are optimal constants
satisfying $\mrm{cond}(\mrm{Q}_{\textsc{nn}} \mR^*\mT_{\mL}\mR) = \lambda_{+}/\lambda_{-}$.
\hfill $\Box$

\quad\\
We end this contribution by presenting results from an actual
numerical experiment testing this preconditioner. We come back to the
concrete setup considered in Section
\ref{NumericIllustrationSec}. This time we solve the preconditioned
linear system $\mrm{Q}_{\textsc{hx}}(\underline{\mR}^*\mT_{\mrm{M}}
\underline{\mR})\bu = \mrm{Q}_{\textsc{hx}}(\bff)$, where $\bff$ is
generated randomly in a similar manner. We simply consider $\gamma =
1$ and, as regards the preconditioner for $\mathscr{L}$, take $\alpha
= \beta = 1$.

\noindent 
\begin{minipage}{0.45\linewidth}  
  \centering\includegraphics[height=6cm]{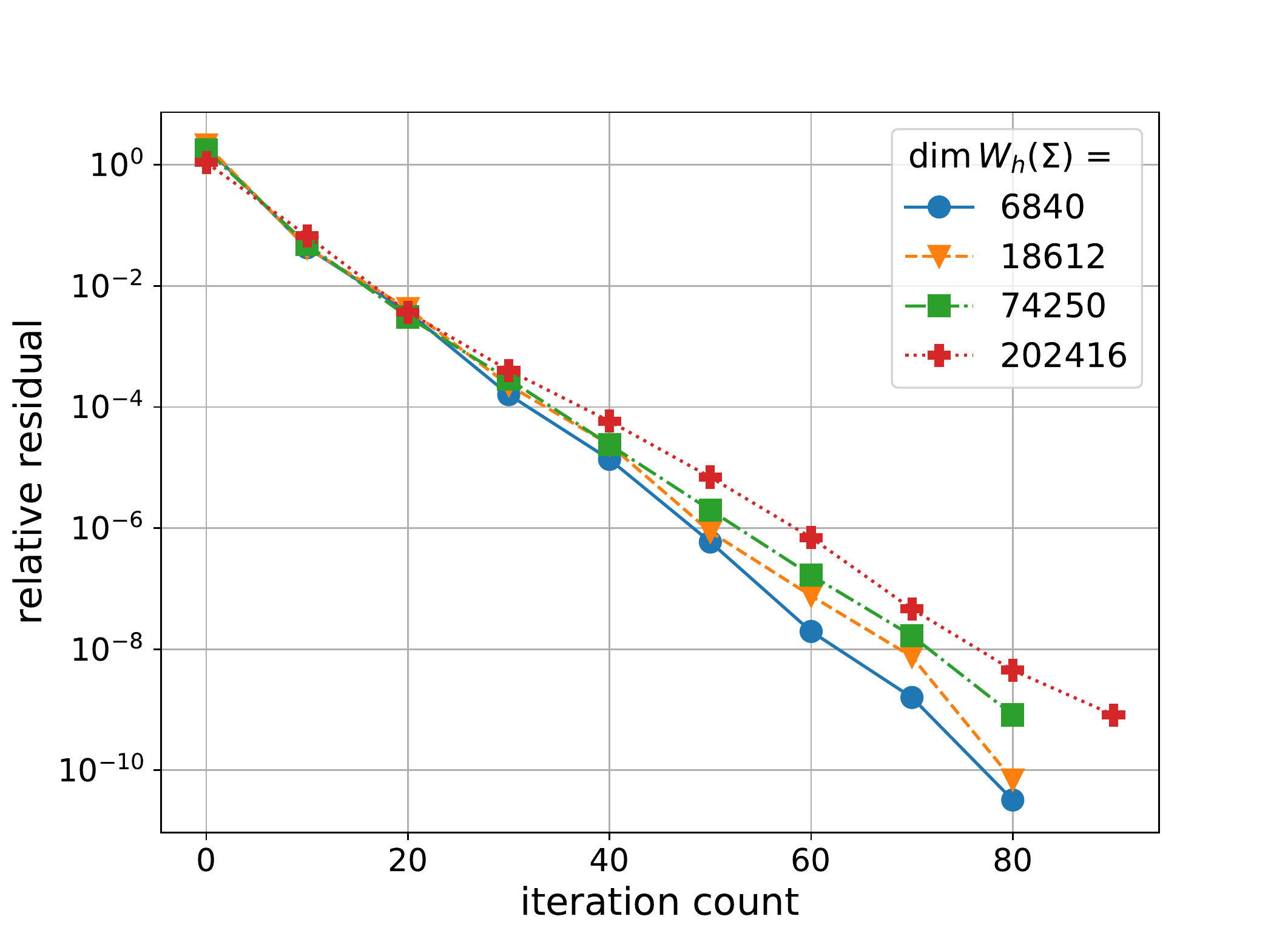}
\end{minipage}
\hspace{2cm}
\begin{minipage}{0.4\linewidth}
  \centering 
  \begin{tabular}{r|r|p{1cm}}
    $\mrm{dim}\,\Wh(\Sigma)$ & $\mrm{dim}\,\Wh(\Omega)$ & \# iter\\
    \hline
    6840   \hspace{0.1cm} &  13915   \hspace{0.1cm} &\hspace{0.3cm} 81\\
    18612  \hspace{0.1cm} &  54472   \hspace{0.1cm} &\hspace{0.3cm} 81\\
    74250  \hspace{0.1cm} &  413388  \hspace{0.1cm} &\hspace{0.3cm} 90\\
    202416 \hspace{0.1cm} &  1796823 \hspace{0.1cm} &\hspace{0.3cm} 97
  \end{tabular}\quad\\[5pt]
  \captionof{table}{Required iterations to reach $10^{-9}$ relative residual
  in PCG applied to $\mR^*\mT_{\textsc{M}}\mR$ preconditioned
  with $\mrm{Q}_{\textsc{hx}}$.} \label{tab:title3} 
  
\end{minipage}

\quad\\[10pt]
For a fixed residual threshold of PCG, the iteration count only
suffers from a mild deterioration that is probably inherited from
$\mQ_{\textsc{nn}}$, in agreement with Lemma \ref{FinalCondEstimate}.


\begin{thebibliography}{10}

\bibitem{MR3649422}
B.~Ayuso~de Dios, R.~Hiptmair, and C.~Pagliantini.
\newblock Auxiliary space preconditioners for {SIP}-{DG} discretizations of
  {${\rm H(curl)}$}-elliptic problems with discontinuous coefficients.
\newblock {\em IMA J. Numer. Anal.}, 37(2):646--686, 2017.

\bibitem{MR4211733}
A.T. Barker and T.~Kolev.
\newblock Matrix-free preconditioning for high-order {$H(\rm curl)$}
  discretizations.
\newblock {\em Numer. Linear Algebra Appl.}, 28(2):Paper No. e2348, 17, 2021.

\bibitem{zbMATH03189639}
A.~Ben-Israel and A.~Charnes.
\newblock Contributions to the theory of generalized inverses.
\newblock {\em J. Soc. Ind. Appl. Math.}, 11:667--699, 1963.

\bibitem{zbMATH01889793}
A.~Ben-Israel and T.N.~E. Greville.
\newblock {\em Generalized inverses. {Theory} and applications.}, volume~15 of
  {\em CMS Books Math./Ouvrages Math. SMC}.
\newblock New York, NY: Springer, 2nd ed. edition, 2003.

\bibitem{MR1192782}
M.~Sh. Birman and M.~Z. Solomjak.
\newblock {\em Spectral theory of selfadjoint operators in {H}ilbert space}.
\newblock Mathematics and its Applications. D. Reidel Publishing Co.,
  Dordrecht, 1987.

\bibitem{MR2970732}
C.R. Dohrmann and O.B. Widlund.
\newblock An iterative substructuring algorithm for two-dimensional problems in
  {$H({\rm curl})$}.
\newblock {\em SIAM J. Numer. Anal.}, 50(3):1004--1028, 2012.

\bibitem{MR3242973}
C.R. Dohrmann and O.B. Widlund.
\newblock Some recent tools and a {BDDC} algorithm for 3{D} problems in
  {$H({\rm curl})$}.
\newblock In {\em Domain decomposition methods in science and engineering
  {XX}}, volume~91 of {\em Lect. Notes Comput. Sci. Eng.}, pages 15--25.
  Springer, Heidelberg, 2013.

\bibitem{MR3465088}
C.R. Dohrmann and O.B. Widlund.
\newblock A {BDDC} algorithm with deluxe scaling for three-dimensional {$H({\bf
  curl})$} problems.
\newblock {\em Comm. Pure Appl. Math.}, 69(4):745--770, 2016.

\bibitem{MR2566786}
C.~Geuzaine and J.-F. Remacle.
\newblock Gmsh: {A} 3-{D} finite element mesh generator with built-in pre- and
  post-processing facilities.
\newblock {\em Internat. J. Numer. Methods Engrg.}, 79(11):1309--1331, 2009.

\bibitem{MR2009375}
R.~Hiptmair.
\newblock Finite elements in computational electromagnetism.
\newblock {\em Acta Numer.}, 11:237--339, 2002.

\bibitem{MR4143285}
R.~Hiptmair and C.~Pechstein.
\newblock A review of regular decompositions of vector fields: continuous,
  discrete, and structure-preserving.
\newblock In {\em Spectral and high order methods for partial differential
  equations---{ICOSAHOM} 2018}, volume 134 of {\em Lect. Notes Comput. Sci.
  Eng.}, pages 45--60. Springer, Cham, 2020.

\bibitem{MR2361899}
R.~Hiptmair and J.~Xu.
\newblock Nodal auxiliary space preconditioning in {${\bf H}({\bf curl})$} and
  {${\bf H}({\rm div})$} spaces.
\newblock {\em SIAM J. Numer. Anal.}, 45(6):2483--2509, 2007.

\bibitem{MR4319101}
Q.~Hu.
\newblock Convergence of the {H}iptmair-{X}u preconditioner for {M}axwell's
  equations with jump coefficients ({I}): extensions of the regular
  decomposition.
\newblock {\em SIAM J. Numer. Anal.}, 59(5):2500--2535, 2021.

\bibitem{MR2536904}
T.V. Kolev and P.S. Vassilevski.
\newblock Parallel auxiliary space {AMG} for {$H({\rm curl})$} problems.
\newblock {\em J. Comput. Math.}, 27(5):604--623, 2009.

\bibitem{li2023nodal}
Y.~Li.
\newblock Nodal auxiliary space preconditioning for the surface de rham
  complex.
\newblock {\em Found. Comput. Math.}, pages 1--30, 2023.

\bibitem{MR2445659}
T.P.A. Mathew.
\newblock {\em Domain decomposition methods for the numerical solution of
  partial differential equations}, volume~61 of {\em Lecture Notes in
  Computational Science and Engineering}.
\newblock Springer-Verlag, Berlin, 2008.

\bibitem{MR2059447}
P.~Monk.
\newblock {\em Finite element methods for {M}axwell's equations}.
\newblock Numerical Mathematics and Scientific Computation. Oxford University
  Press, New York, 2003.

\bibitem{MR3013465}
C.~Pechstein.
\newblock {\em Finite and boundary element tearing and interconnecting solvers
  for multiscale problems}, volume~90 of {\em Lecture Notes in Computational
  Science and Engineering}.
\newblock Springer, Heidelberg, 2013.

\bibitem{zbMATH03302838}
W.~V. Petryshyn.
\newblock On generalized inverses and on the uniform convergence of
  {{\((I-\beta K)^n\)}} with application to iterative methods.
\newblock {\em J. Math. Anal. Appl.}, 18:417--439, 1967.

\bibitem{zbMATH01022519}
W.~Rudin.
\newblock {\em Functional analysis.}
\newblock New York, NY: McGraw-Hill, 2nd ed. edition, 1991.

\bibitem{zbMATH01953446}
Y.~Saad.
\newblock {\em Iterative methods for sparse linear systems}.
\newblock Boston, MA: PWS Publishing Company, 1996.

\bibitem{MR1410757}
B.F. Smith, P.E. Bjorstad, and W.D. Gropp.
\newblock {\em Domain decomposition}.
\newblock Cambridge University Press, Cambridge, 1996.
\newblock Parallel multilevel methods for elliptic partial differential
  equations.

\bibitem{MR2104179}
A.~Toselli and O.~Widlund.
\newblock {\em Domain decomposition methods---algorithms and theory}, volume~34
  of {\em Springer Series in Computational Mathematics}.
\newblock Springer-Verlag, Berlin, 2005.

\bibitem{MR3643544}
O.B. Widlund and C.R. Dohrmann.
\newblock B{DDC} deluxe domain decomposition.
\newblock In {\em Domain decomposition methods in science and engineering
  {XXII}}, volume 104 of {\em Lect. Notes Comput. Sci. Eng.}, pages 93--103.
  Springer, Cham, 2016.

\bibitem{MR1393008}
J.~Xu.
\newblock The auxiliary space method and optimal multigrid preconditioning
  techniques for unstructured grids.
\newblock volume~56, pages 215--235. 1996.
\newblock International GAMM-Workshop on Multi-level Methods (Meisdorf, 1994).

\end{thebibliography}

\end{document}